\documentclass[11pt]{amsart}
\usepackage{amssymb,mathrsfs,graphicx,enumerate,color}
\usepackage{amsmath, amsfonts, amscd, epsfig}

\usepackage{mathtools}
\usepackage{hyperref}
\usepackage{comment}

\newcommand{\Ccal}{\mathcal{C}}

\newcommand{\Ocal}{\mathcal{O}}

\renewcommand{\hbar}{\bar{h}}

\newcommand{\ptil}{\tilde{p}}

\newcommand{\stil}{\tilde{s}}

\newcommand{\util}{\tilde{u}}
\newcommand{\vtil}{\tilde{v}}

\newcommand{\abs}[1]{\left|#1\right|}

\DeclareMathOperator{\tr}{tr}

\renewcommand{\b}{\beta}
\newcommand{\g}{\gamma}

\newcommand{\e}{\varepsilon}

\renewcommand{\th}{\theta}
\renewcommand{\k}{\kappa}

\newcommand{\m}{\mu}

\newcommand{\x}{\xi}

\newcommand{\s}{\sigma}
\renewcommand{\t}{\tau}

\newcommand{\rd}{\partial}

\newcommand{\thtil}{\tilde{\th}}

\numberwithin{equation}{section}

\usepackage{colortbl}
\definecolor{black}{rgb}{0.0, 0.0, 0.0}
\definecolor{red}{rgb}{1.0, 0.5, 0.5}

\title[   ]{Traveling Wave Solutions to a Large Class of Brenner-Navier-Stokes-Fourier Systems}

\author[Eo]{Saehoon Eo}
\address[Saehoon Eo]
{ Department of Mathematics, \newline
Stanford University \\ CA 94305, USA}
\email{eosehoon@stanford.edu}

\author[Eun]{Namhyun Eun}
\address[Namhyun Eun]
{ Department of Mathematical Sciences, \newline
Korea Advanced Institute of
Science and Technology \\
Daejeon 34141, Korea}
\email{namhyuneun@kaist.ac.kr}
 
\newtheorem{theorem}{Theorem}[section]

\newtheorem{corollary}{Corollary}[section]
\newtheorem{proposition}{Proposition}[section]
\newtheorem{remark}{Remark}[section]

\newtheorem{definition}{Definition}[section]

\newcommand{\RR}{{\mathbb{R}}}

\newcommand{\vt}{{\tilde{v}_\e}}
\newcommand{\ut}{{\tilde{u}_\e}}

\newcommand{\tht}{{\tilde{\th}_\e}}

\newcommand{\step}[1]{\vskip0.2cm \noindent{\it Step #1:} }

\begin{document}

\date{\today}


\thanks{\textbf{Acknowledgment.}
The authors thank Professor Moon-Jin Kang for his valuable comments and continuous support throughout the preparation of this work.
This work was supported by Samsung Science and Technology Foundation under Project Number SSTF-BA2102-01. }

\begin{abstract}
The Brenner-Navier-Stokes-Fourier (BNSF) system, introduced by Howard Brenner, was developed to address some deficiencies in the classical Navier-Stokes-Fourier system, based on the concept of volume velocity.
We consider the one-dimensional BNSF system in Lagrangian mass coordinates, incorporating temperature-dependent transport coefficients, which yields a more physically realistic framework.
We establish the existence and uniqueness of monotone traveling wave solutions (or viscous shocks) to the BNSF system with any positive \(\Ccal^2\) dissipation coefficients, provided that the shock amplitude is sufficiently small.
We utilize geometric singular perturbation theory as in the constant coefficient case \cite{EEKO}; however, due to the arbitrary nonlinearities of the coefficients, we employ the implicit function theorem, which grants robustness to our approach.
This work is motivated by \cite{EEK}, which proves a contraction property of any large solutions to the BNSF system around the traveling wave solutions.
Thus, we also derive some quantitative estimates on the traveling wave solutions that play a fundamental role in \cite{EEK}.
\end{abstract}

\maketitle \centerline{\date}

\tableofcontents

\section{Introduction}
\setcounter{equation}{0}
We consider the Brenner-Navier-Stokes-Fourier(BNSF) system with one spatial variable. 
In the Lagrangian mass coordinates, it takes the following form (see \cite{EEK}):
\begin{equation}
\begin{cases} \label{main0}
v_t-(u_m)_x = 0, \\
(u_v)_t+p(v,\th)_x = \Big(\m(\th)\frac{(u_v)_x}{v}\Big)_x, \\
E_t+(p(v,\th)u_v)_x = 
\Big(\k(\th)\frac{\th_x}{v}\Big)_x+\Big(\m(\th)\frac{u_v(u_v)_x}{v}\Big)_x.
\end{cases}
\end{equation}
Here, \(v\) denotes the specific volume, \(u_m\) mass velocity, \(u_v\) volume velocity, \(\th\) absolute temperature, \(p\) pressure law, and \(E\) the total energy.
The concepts of mass and volume velocity are described in detail below.
Moreover, \(\m\) and \(\k\) denote the viscosity and the heat conductivity, respectively, both of which are positive functions of absolute temperature \(\th\).

\vspace{2mm}
The BNSF system has been proposed by Brenner \cite{brenner2006fluid} to improve some flaws of the compressible Navier-Stokes-Fourier(NSF) system.
He mentioned in \cite{brenner2005kinematics,brenner2005navier,brenner2006fluid} that the NSF system does not adequately describe compressible flows, particularly under extreme conditions such as rarefied gases, shock waves, and gaseous flows through micro-channels.
Indeed, a substantial number of studies have identified deficiencies in the NSF system and proposed various improvements (see \cite{arkilic2001mass, chakraborty2007derivations, dadzie2008continuum, dongari2010predicting, dongari2009extended, durst2011treatments, greenshields2007structure, harley1995gas, klimontovich1992need, klimontovich1993hamiltonian, reese2003new}), but Brenner's justification is regarded as the most systematic and well-established account.

\vspace{2mm}
His idea is founded on the bi-velocity theory, which claims the existence of two distinct velocities when describing fluid dynamics, namely, the mass velocity \(u_m\) (which is identical to the fluid velocity in classical equations) and the volume velocity \(u_v\).
In \cite{brenner2006fluid}, he contended that these two velocities are implicitly assumed to be identical in the bulk of the literature on continuum fluids, and that this assumption has largely remained unquestioned.
However, he observed that, in general, \(u_m \neq u_v\) and that the discrepancy becomes more significant as the density (or the specific volume) gradient increases.
More explicitly, in the Lagrangian variables, the constitutive relation between \(u_v\) and \(u_m\) suggested by Brenner \cite{brenner2006fluid} is as follows:
\begin{equation} \label{umuv}
u_m=u_v+\frac{\k(\th)}{c_p}\Big(\frac{v_x}{v}\Big)
\end{equation}
where \(c_p\) is the specific heat at constant pressure.
The difference between the two velocities is proportional to \(v_x\), and thus it is worth mentioning that the concept of volume velocity becomes insignificant in incompressible fluids with uniform density.

\vspace{2mm}
In a series of works \cite{brenner2005kinematics,brenner2005navier,brenner2006fluid}, Brenner proposed the constitutive relation \eqref{umuv} based on the theory by \"Ottinger and justified it through a variety of ways including the alignment between Burnett's solution to the Boltzmann equation and empirical data.
He also established in \cite{brenner2006fluid} that the relation \eqref{umuv} is consistent with linear irreversible thermodynamics, and remarkably, for the case of fluids in which density gradients exist, no-slip boundary condition should be imposed in terms of volume velocity to achieve better agreement with experimental data.

\vspace{2mm}
The motivation for introducing volume velocity can be understood from a more intuitive perspective as follows.
When gas particles are displaced due to internal pressure, work is done.
To quantify this work, one must define the displacement distance.
However, it is not individual gas particles that generate pressure and are displaced by it, but rather a collection of particles.
Since the boundaries of such collections cannot be clearly defined, it intuitively appears impossible to rigorously define the displacement distance.
This ambiguity reflects the fact that the volume of gas is not a point-based property; rather, it is a statistical quantity associated with a collection of particles.
Consequently, the concept of volume velocity, as a new form of velocity, is necessary to describe the motion of particle collections, rather than relying on the traditional concept of velocity.
To put it briefly, the mass velocity \(u_m\) is used for the mass transportation and convective phenomena, while the volume velocity is introduced to describe the momentum, energy, work and viscous stress.

\vspace{2mm}
While Brenner's theory is thorough, mathematical studies on the BNSF system are still scarce.
Feireisl-Vasseur \cite{FeVa} established the global existence of weak solutions for the multi-dimensional initial boundary value problem.
The BNSF system is also used to construct measure-valued solutions of Euler in \cite{Fe1}, and the finite volume scheme for Euler in \cite{Fe2}.

\vspace{2mm}
Thanks to the constitutive relation \eqref{umuv}, the system \eqref{main0} can be written as
\begin{equation}
\left\{
\begin{aligned} \label{main}
    &v_t-u_x = \Big(\t(\th)\frac{v_x}{v}\Big)_x, \\
    &u_t+p(v,\th)_x = \Big(\m(\th)\frac{u_x}{v}\Big)_x, \\
    &E_t+(p(v,\th)u)_x = 
    \Big(\k(\th)\frac{\th_x}{v}\Big)_x+\Big(\m(\th)\frac{uu_x}{v}\Big)_x,
\end{aligned} \right.
\end{equation}
where \(u\) is now the volume velocity and \(\t(\th) = \k(\th)/c_p\) is the Brenner coefficient.
The total energy is given as \(E=e+\frac{u^2}{2}\), where \(e\) is the internal energy.
We also consider the ideal polytropic gas where the pressure \(p\) and the internal energy \(e\) are given by
\begin{equation} \label{pressure}
p(v,\th)=\frac{R\th}{v}, \qquad e(\th)=\frac{R}{\g-1}\th,
\end{equation}
where \(R>0\) is the gas constant and \(\g>1\) is the adiabatic constant.

The main advancement of the present paper over the previous work \cite{EEKO} lies in the incorporation of temperature-dependent transport coefficients, i.e., viscosity and heat conductivity.
It is well-known from the Chapman-Enskog theory or the first level of approximation in kinetic theory that the viscosity and the heat conductivity coefficients are functions of temperature alone.
According to Chapman-Cowling \cite{ChapmanCowling90} or Vincenti-Kruger \cite{VincentiKruger66}, for a gas with intermolecular potential proportional to \(r^{-a}\), where \(r\) is the intermolecular distance, \(\m(\th)\) and \(\k(\th)\) are given as follows: for some positive constants \(\m_0\) and \(\k_0\),
\begin{equation*}
\m(\th) = \m_0\th^\b, \qquad \k(\th) = \k_0\th^\b, \qquad \text{ where } \b = \frac{a+4}{2a}.
\end{equation*}
Thus, the Brenner coefficient \(\t(\th) = \k(\th)/c_p\) may also degenerate as \(\theta^\beta\) near zero temperature.
However, the situation is more subtle in practice.
While the specific heat \(c_p\) is often treated as a constant, it is in fact known to degenerate as \(\th\to0\) (see \cite{kittel1980thermal,kittel2018introduction}).
Therefore, we take into account arbitrary positive functions \(\t(\th),\m(\th)\) and \(\k(\th)\) with \(\Ccal^2\)-regularity, and prove the existence and uniqueness of the traveling wave solutions.

\vspace{2mm}
This work is motivated by \cite{EEK}, which resolves a challenging open problem in the field of hyperbolic conservation laws concerning the stability of irreversible singularities, in particular, Riemann shocks, within the class inviscid limits from associated physical viscous systems (see \cite{KV-Inven}). Eo-Eun-Kang \cite{EEK} established contraction properties for arbitrarily large solutions to the BNSF system with \(\t(\th)=1+\th^2,\m(\th)=\th^2,\k(\th)=\th^2\), perturbed from a viscous shock.
From the contraction estimates, the existence of inviscid limits and the stability of the corresponding Riemann shock could be obtained.
The existence of a traveling wave solution---also known as a viscous shock---and its quantitative estimates play a central role in \cite{EEK}.
For the NSF system, we refer to \cite{LYZZ} for the consideration of temperature-dependent viscosity and heat conductivity.

\vspace{2mm}
The approach proposed in the present paper is expected to be robust.
Unlike the constant dissipation case \cite{EEKO}, where it is possible to eliminate the \(\th\) variable and reformulate the system solely in terms of the \(v\) variable, the general case does not allow such a simplification.
To address this issue and to uncover the necessary geometric structure of the system, we make an effective use of the implicit function theorem, which has a broad applicability in nonlinear analysis.
This will be discussed in Section \ref{subsec:idea}.
Regarding the applicability, it was shown in \cite{HKKL} that the Navier-Stokes-Korteweg system with constant viscosity admits a traveling wave solution.
It is therefore expected that our approach can extend it to the case of non-constant viscosity, and even to the Navier-Stokes-Fourier-Korteweg system \cite{Plaza} with temperature-dependent transport coefficients.
Moreover, due to the lack of nonlinearity in the continuity equation, our approach is also anticipated to be applicable for the existence of traveling wave solutions for the NSF system with temperature-dependent coefficients.

\subsection{Main results}
In this paper, we establish that the system \eqref{main} admits tor viscous shocks which connect two end states \((v_-,u_-,E_-)\) and \((v_+,u_+,E_+)\) when these two states are close enough and satisfy the Rankine-Hugoniot condition and Lax entropy condition as follows:
\begin{equation}
\begin{aligned} \label{end-con}
&\exists~\s \quad \text{s.t.}~\left\{
\begin{aligned}
&-\s(v_+-v_-) -(u_+-u_-) =0, \\
&-\s(u_+-u_-) +p(v_+,\th_+)-p(v_-,\th_-) =0, \\
&-\s(E_+-E_-) +p(v_+,\th_+)u_+ -p(v_-,\th_-)u_- =0,
\end{aligned} \right. \\
&\text{and either \(v_->v_+\), \(u_->u_+\), \(\th_-<\th_+\) or \(v_-<v_+\), \(u_->u_+\), \(\th_->\th_+\) holds.}
\end{aligned}
\end{equation}
Here, \(\th_-\) and \(\th_+\) satisfy \(E_- = \frac{R}{\g-1}\th_- + \frac{u_-^2}{2}\) and \(E_+ = \frac{R}{\g-1}\th_+ + \frac{u_+^2}{2}\).
More precisely, for given constant states \((v_-,u_-,E_-)\) and \((v_+,u_+,E_+)\) satisfying \eqref{end-con}, provided that \(\abs{v_+-v_-}\) is sufficiently small, we prove the existence of a viscous shock wave (\(\vtil,\util,\thtil)(\x)=(\vtil, \util, \thtil\))(\(x-\s t\)), as a traveling wave solution to the following system of ODEs:
\begin{equation}
\begin{cases} \label{shock_0}
    -\s \vtil'-\util' = \Big(\t(\thtil)\frac{\vtil'}{\vtil}\Big)', \\
    -\s \util'+p(\vtil,\thtil)' = \Big(\m(\thtil)\frac{\util'}{\vtil}\Big)', \\
    -\s \Big(\frac{R}{\g-1}\thtil + \frac{\util^2}{2}\Big)'+(p(\vtil,\thtil)\util)' = 
    \Big(\k(\thtil)\frac{\thtil'}{\vtil}\Big)'+ \Big(\m(\thtil)\frac{\util \util'}{\vtil}\Big)', \\
    \lim_{\x\to\pm\infty} (\vtil,\util,\thtil)(\x)=(v_\pm,u_\pm,\th_\pm).
\end{cases}
\qquad\Big(\,\,'=\frac{d}{d\x}\Big)
\end{equation}
Here, if \(v_->v_+\), (\(\vtil, \util, \thtil\))(\(x-\s t\)) is a 1-shock wave with velocity \(\s=- \sqrt{-\frac{p_+-p_-}{v_+-v_-}}<0\), where \(p_\pm := p(v_\pm, \th_\pm)\).
If \(v_-<v_+\), it is a 3-shock wave with \(\s= \sqrt{-\frac{p_+-p_-}{v_+-v_-}}>0\).

\vspace{2mm}
Our first result is on the existence and uniqueness of traveling wave solutions to \eqref{main} as monotone profiles satisfying \eqref{shock_0}, which is stated in the following theorem.

\begin{theorem}\label{thm_EU} (Existence and Uniqueness)
Let \(\t(\th),\m(\th),\k(\th)\) be any positive \(\Ccal^2\)-functions.
For a given left-end state \((v_-,u_-,\th_-) \in \RR^+ \times \RR \times \RR^+\), there exists a constant \(\e_0>0\) such that for any right-end state \((v_+, u_+, \th_+) \in \RR^+ \times \RR \times \RR^+\) satisfying \eqref{end-con} and \(\e := \abs{v_--v_+} \in (0,\e_0)\), there is a unique traveling wave solution \((\vt, \ut, \tht)\colon \RR\to \RR^+ \times \RR \times \RR^+\) to \eqref{main}, as a monotone profile satisfying \eqref{shock_0}, which satisfies \eqref{1shock} when \(v_->v_+\) or \eqref{derivsign} when \(v_-<v_+\).
\end{theorem}

The next theorem provides certain quantitative estimates of traveling wave solutions. 
In particular, we present the estimates of the ratio among \(\vtil',\util',\thtil'\), and the exponential decay of the waves.
This is essential for the uniform stability estimates for traveling waves \cite{EEK}.

\begin{theorem}\label{thm_estimates} (Quantitative Estimates)
Let \(\t(\th),\m(\th),\k(\th)\) be any positive \(\Ccal^2\)-functions.
For a given left-end state \((v_-,u_-,\th_-) \in \RR^+ \times \RR \times \RR^+\), there exist constants \(\e_0,C,C_1>0\) such that for any right-end state \((v_+, u_+, \th_+) \in \RR^+ \times \RR \times \RR^+\) with \eqref{end-con}, \(\e := \abs{v_--v_+} \in (0,\e_0)\), the following holds.
Let \((\vt, \ut, \tht)\) be the monotone solution to \eqref{shock_0} with \(\vt(0) = \frac{v_-+v_+}{2}\).\\
Then, the solution \((\vt, \ut, \tht)\) satisfies the following estimates:
\begin{align}
&\begin{aligned}\label{decay}
&\abs{(\vt(\x)-v_-, \ut(\x)-u_-, \tht(\x)-\th_-)} \le C\e e^{-C_1\e\abs{\x}} \text{ for all } \x \le 0, \\
&\abs{(\vt(\x)-v_+, \ut(\x)-u_+, \tht(\x)-\th_+)} \le C\e e^{-C_1\e\abs{\x}} \text{ for all } \x \ge 0, \\
\end{aligned} \\
&\abs{(\vtil'_\e(\x), \util'_\e(\x), \thtil'_\e(\x))} \le C\e^2e^{-C_1\e\abs{\x}} \text{ for all } \x\in \RR, \label{derivdecay}\\
&\abs{(\vtil''_\e(\x), \util''_\e(\x), \thtil''_\e(\x))} \le C\e\abs{(\vtil'_\e(\x), \util'_\e(\x), \thtil'_\e(\x))} \text{ for all } \x\in \RR.
\label{2bound1}
\end{align}
It also holds that \(\abs{\vt'} \sim \abs{\ut'} \sim |\tht'|\) for all \(\x\in\RR\), and more explicitly, 
\begin{align}
\abs{\util'_\e(\x) + \s_*\vtil'_\e(\x)} &\le C\e\abs{\vtil'_\e(\x)}, 
\label{ratio-vu} \\
\abs{\thtil'_\e(\x) + \frac{(\g-1)p_-}{R}\vtil'_\e(\x)} &\le C\e\abs{\vtil'_\e(\x)},
\label{ratio-vth}
\end{align}
where
\[
\s_\e := -\sqrt{-\frac{p_+-p_-}{v_+-v_-}}   \quad\text{ and }\quad \s_* := -\sqrt{\frac{\g p_-}{v_-}} = -\frac{\sqrt{\g R \th_-}}{v_-},
\]
or
\[
\s_\e := \sqrt{-\frac{p_+-p_-}{v_+-v_-}}  \quad\text{ and }\quad \s_* := \sqrt{\frac{\g p_-}{v_-}} = \frac{\sqrt{\g R \th_-}}{v_-},
\]
which satisfy
\begin{equation} \label{sm1}
\abs{\s_\e-\s_*} \le C \e.
\end{equation}
In addition, if we consider the 1-shock, i.e., \(\s_\e<0\), then we have
\begin{equation} \label{1shock}
\vtil'_\e <0, \quad \util'_\e <0, \quad \thtil'_\e >0 \text{ for all } \x \in \RR,
\end{equation}
and for the 3-shock case, i.e., \(\s_\e>0\), we have
\begin{equation}\label{derivsign}
\vtil'_\e >0, \quad \util'_\e <0, \quad \thtil'_\e <0 \text{ for all } \x \in \RR.
\end{equation}
Furthermore, we define a constant \(A>0\) by
\begin{equation}\label{Adef}
A \coloneqq \frac{R\g\s_*}{R\t(\th_-)+R\g\m(\th_-)+(\g-1)^2\k(\th_-)}\frac{\g+1}{2},
\end{equation}
and it satisfies the following estimates holds: 
\begin{equation} \label{tail}
\abs{\vt'- A(v_+-\vt)(\vt-v_-)} \le C\e(v_+-\vt)(\vt-v_-).
\end{equation}
\end{theorem}

The last estimate \eqref{tail} implies the following:
\begin{corollary}\label{Jaccor}
For \(y(\x) \coloneqq (\vt(\x)-v_-)/\e\), it follows that
\begin{equation}\label{Jacineq}
\abs{\frac{1}{y(1-y)}\frac{dy}{d\x} - A\e} \le C\e^2.
\end{equation}
\end{corollary}

We conclude this section with a few remarks.
\begin{remark}
Theorem \ref{thm_EU} and \ref{thm_estimates} significantly generalizes the previous work \cite{EEKO} for the constant dissipation coefficient case.
Within the present framework, additional challenges arise, as the argument in \cite{EEKO} crucially relies on the fact that the continuity equation does not involve the temperature variable.
A detailed comparison between the present work and \cite{EEKO} is provided in Section \ref{sec:pre}.
\end{remark}

\begin{remark}
The set of estimates in Theorem \ref{thm_estimates} is essential for the stability estimates of traveling wave solutions, based on the method of a-contraction with shifts \cite{KV21}, as in the previous studies \cite{HKK23,HKKL,KV-Inven,KV-2shock,KVW23,KVW-NSF}.
To be precise, we refer to \cite[Lemma 2.1]{KV21} and \cite[Lemma 2.3]{KVW-NSF}.
Recently, Eo-Eun-Kang \cite{EEK} established a certain contraction property of large solutions to the BNSF system with temperature-dependent coefficient perturbed from viscous shocks, which requires the estimates in Theorem \ref{thm_estimates} in this exact form.
\end{remark}

\begin{remark}
In contrast to the previous work \cite{EEKO} of constant coefficient case, we determine the constant \(A\) that characterizes the ratio between two quantities \(\vt'\) and \((\vt-v_-)(v_+-\vt)\), up to an \(\Ocal(\e)\) error.
It is obvious by \eqref{tail} that the weaker inequality holds:
\begin{equation*}
C^{-1}(\vt-v_-)(v_+-\vt) \le \vt'\le C(\vt-v_-)(v_+-\vt),
\end{equation*}
which is essential to establish \cite[Lemma 4.2 and 4.3]{EEK}, which were the key lemmas for the contraction property.
Moreover, \eqref{tail} also yields Corollary \ref{Jaccor}, which coincides with that of \cite[Lemma 5.1]{EEK}, and the proof is given at the end of Section \ref{sec:property}.
Indeed, the constant \(A\) serves as a coefficient associated with the critical manifold \eqref{Critmfd}.
\end{remark}

\begin{remark}
It suffices to prove Theorem \ref{thm_estimates} for 3-shock cases.
In fact, the result for 1-shocks can be obtained by the change of variables \(x \mapsto -x, u \mapsto -u, \s_\e \mapsto -\s_\e\) and \(\s_* \mapsto -\s_*\). Therefore, in the sequel, we only consider the 3-shock case, i.e.,
\[
\s_\e = \sqrt{-\frac{p_+-p_-}{v_+-v_-}}>0 \quad\text{ and }\quad
\s_* \coloneqq \sqrt{\frac{\g p_-}{v_-}} = \frac{\sqrt{\g R \th_-}}{v_-}.
\]    
\end{remark}

\section{Key ideas and methodology} \label{sec:pre}
\setcounter{equation}{0}

\subsection{Key idea behind the proofs} \label{subsec:idea}
First of all, we integrate the ODE system \eqref{shock_0} over \((-\infty,\x]\) and then multiplying the second equation by \(\util\) and subtracting it from the third equation, we obtain the system of ODEs of first order:
\begin{equation}
\label{intODE}
    \left\{
    \begin{aligned}
        & -\s_\e(\vtil-v_-) - (\util-u_-) = \t(\thtil) \frac{\vtil'}{\vtil}, \\
        & -\s_\e(\util-u_-) + (\ptil-p_-) = \m(\thtil) \frac{\util'}{\vtil}, \\
        &  -\s_\e\Big(\frac{R}{\g-1}(\thtil-\th_-) - \frac{1}{2}(\util-u_-)^2\Big)
    + p_-(\util-u_-) = \k(\thtil) \frac{\thtil'}{\vtil}.
    \end{aligned}
    \right.
\end{equation}
As in \cite{EEKO}, we will use the theory of geometric singular perturbation given by Fenichel's first theorem to prove the existence of the solutions and their estimates.
To this end, we aim to isolate \(\util\) variable from the first equation.
Since we have
\begin{align}\label{uexplicit}
\util-u_- &= -\s_\e(\vtil-v_-) - \t(\thtil) \frac{\vtil'}{\vtil}, \\
\label{utexplicit}
\util' &= -\s_\e \vtil' - \t'(\thtil)\frac{\thtil'\vtil'}{\vtil} - \t(\thtil)\frac{\vtil''\vtil-(\vtil')^2}{\vtil^2},
\end{align}
the system \eqref{intODE} can be reduced into the following system with two variables:
\begin{equation} \label{vthODE}
\left\{
\begin{aligned}
    & 0 = \s_\e\Big(\s_\e(\vtil-v_-) +\t(\thtil)\frac{\vtil'}{\vtil}\Big) + (\ptil - p_-)
    + \frac{\m(\thtil)}{\vtil} \Big(\s_\e \vtil' + \t'(\thtil)\frac{\thtil'\vtil'}{\vtil} + \t(\thtil)\frac{\vtil''\vtil-(\vtil')^2}{\vtil^2}\Big) \\
    &\qquad \eqqcolon h_1(\vtil, \thtil, \vtil', \vtil'', \thtil, \e), \\
    &0 = -\frac{\s_\e R}{\g-1}(\thtil-\th_-) +\frac{\s_\e}{2}\Big(\s_\e(\vtil-v_-)+ \t(\thtil) \frac{\vtil'}{\vtil}\Big)^2
    - p_-\Big(\s_\e(\vtil-v_-) + \t(\thtil) \frac{\vtil'}{\vtil}\Big) - \k(\thtil) \frac{\thtil'}{\vtil} \\
    &\qquad \eqqcolon h_2(\vtil, \thtil, \vtil', \vtil'', \thtil, \e).
\end{aligned}
\right.
\end{equation}
Note that since \(\s_\e\) is a smooth function of \(\e\) near \(0\) (see \cite[(3.17)]{EEKO}) as
\begin{equation}\label{seval}
\s_\e = \frac{\sqrt{\g p_-}}{\sqrt{v_-+\frac{\g+1}{2}\e}},
\end{equation}
the two functions \(h_1\) and \(h_2\) are smooth in \(\e\) as well.

To use Fenichel's first theorem, we need to write the system \eqref{vthODE} as a higher order ODE in \(\vtil\) only, as in \cite{EEKO}.
However, such an approach is not applicable here, since \(\thtil\) cannot be isolated from the first equation.
Therefore, we rather apply the implicit function theorem, to express \(\thtil\) and \(\thtil'\) in terms of \((\vtil, \vtil, \vtil'', \e)\),
i.e., regarding \(h_1\) and \(h_2\) to be functions on \(\RR^6\), and find \(g_1, g_2\colon U_0 \subset \RR^4\to \RR\) which satisfy 
\begin{equation}\label{gdef}
\begin{aligned}
&h_1(\vtil, g_1(\vtil, \vtil', \vtil'', \e), \vtil', g_2(\vtil, \vtil', \vtil'', \e), \vtil'', \e) = 0, \\
&h_2(\vtil, g_1(\vtil, \vtil', \vtil'', \e), \vtil', g_2(\vtil, \vtil', \vtil'', \e), \vtil'', \e) = 0,
\end{aligned}
\end{equation}
where \(U_0\) is a neighborhood of \((v_-,0,0,0)\).
Then, using 
\[
\frac{d}{d\x} g_1(\vtil, \vtil', \vtil'', \e)
= g_2(\vtil, \vtil', \vtil'', \e),
\]
we derive an ODE of \(\vtil\) only as follows:
\begin{equation}\label{absODE}
\vtil''' = \Big(\frac{\rd g_1}{\rd \vtil''}\Big)^{-1}\Big[g_2 - \frac{\rd g_1}{\rd\vtil}\vtil' - \frac{\rd g_1}{\rd \vtil'}\vtil''\Big].
\end{equation}
Moreover, as in \cite{EEKO}, we will use the following transformation:
\begin{equation}\label{vtilinw}
\begin{aligned}
\vtil(\x) = v_- + \e w_0(\e\x), \quad
\vtil'(\x) = \e^2w_1(\e\x), \quad
\vtil''(\x) = \e^2w_2(\e\x), \quad
\vtil'''(\x) = \e^3w_2'(\e\x).
\end{aligned}
\end{equation}
We simply obtain \(w_0' = w_1\) and \(\e w_1' = w_2\), and further, \(\e w_2'\) can be expressed as a function of \(w_0, w_1\) and \(w_2\) by \eqref{absODE}.
This yields an ODE system of the form as \eqref{sys}.
Then we apply the geometric singular perturbation theory in the standard way to the system of \((w_0,w_1,w_2)\).
It is worth mentioning that while the explicit representation of the system cannot be obtained, it is sufficient to compute the first and second derivatives of \(g_1\) and \(g_2\) \textbf{at} \(\mathbf{(v_-,0,0,0)}\) in order to characterize the shape of the critical manifold.
Section \ref{sec:EU} is dedicated to detailing the derivation of the system to which Fenichel's first theorem can be applied.

\subsection{Fenichel's first theorem and locally invariant manifolds}
In this subsection, we present the main concepts for the geometric approach.
We consider the following system of ODEs with a small parameter \(\e\): 
\begin{equation}\label{slow}
\left\{
\begin{aligned}
    \e x' &= f(x, y, \e), \\
    y' &= g(x, y, \e), 
\end{aligned}
\right.
\end{equation}
where \(x\in \RR^n\), \(y\in \RR^l\). Here, we assume \(f\) and \(g\) to be smooth on a set \(U\times I\), where \(U\subset \RR^{n+l}\) is open and \(I\) is an open interval which contains \(0\). 
Following \cite{JONES}, we say an \(l\)-dimensional manifold \(M_0\subset \RR^{n+l}\) is a \textit{critical manifold} of the system \eqref{slow} if each \((x, y)\in M_0\) satisfies that \(f(x, y, 0) = 0\).
To state Fenichel's first theorem, we introduce the following definitions:
\begin{definition}\label{def-normally_hyperbolic}
We say that the manifold \(M_0\) is normally hyperbolic relative to \eqref{slow} if the \(n\times n\) matrix 
\[
\left. D_x f(x, y, \e) \right|_{\e=0}
\]
has \(n\) eigenvalues (counting multiplicity) with nonzero real part for each \((x, y)\in M_0\). 
\end{definition}

Fenichel's first theorem establishes the existence of a manifold to which the flow of solutions to \eqref{slow} is confined.
To be precise, we first introduce the notion of a locally invariant manifold:
\begin{definition}\label{def-locally_invariant}
    A set \(M \subseteq \RR^{n+l}\) is locally invariant under the flow from \eqref{slow} if there is an open set \(V \subseteq \RR^{n+l}\) containing \(M\) such that for any \(x_0\in M\) and any \(T\in\RR\), the trajectory \(x\) starting at \(x_0\) with \(x([0, T])\subset V\) also satisfies \(x([0, T])\subset M\). (When \(T<0\), replace \([0,T]\) with \([T,0]\)). 
\end{definition}

Note that, if \(M_0\) is normally hyperbolic, then by the implicit function theorem, it can be locally represented as a graph of \(x\) in terms of \(y\).
However, in the following discussion, we restrict our attention to the case in which \(M_0\) can be globally expressed as a graph of \(x\) in terms of \(y\) over a compact domain \(K\subset \RR^l\).
Moreover, in this case, Fenichel's first theorem can be stated as follows.
\begin{proposition}  \cite{JONES} \label{prop_FFT}
    Consider the system \eqref{slow} with normally hyperbolic critical manifold \(M_0\).
    Suppose that there is a smooth function \(h\colon K\to \RR^n\) whose graph is contained in \(M_0\), i.e.,
    \[
    \{(x, y)\mid x = h(y), \,y\in K\} \subset M_0
    \]
    where \(K\subset \RR^l\) is a compact and simply connected smooth set. 
    Then, there exists \(\e_0>0\) such that for each \(\e\) with \(\abs{\e}<\e_0\), there is a function \(h^\e\colon K\to \RR^n\) which satisfies the followings: for each \(\e \neq 0\), the graph
    \[
    \{(x, y)\mid x = h^\e(y), \,y\in K\}
    \]
    is contained in an \(l\)-dimensional manifold \(M_\e\) which is locally invariant under the flow of \eqref{slow}. Moreover, \(h^\e\) can be taken to be \(\Ccal^r\) for each \(r < \infty\), jointly in \(y\) and \(\e\) and \(h=h^0\) on \(K\) where $h^0$ denotes \(h^{\e}\) at \(\e=0\). 
\end{proposition}
\begin{remark} \label{rmk-FFT}
    Since the family of functions \(h^\e\) is regular in \(\e\) as well, a direct consequence of Proposition \ref{prop_FFT} is that there exists a function \(s\) with two variables \(\e\) and \(y\) which satisfies 
    \[
    h^\e(y) = h(y) + \e s(y, \e).
    \]
    Note that the function \(s\) has the same regularity as \(h^\e\), namely, \(s\) is a \(\Ccal^r\)-function for each \(r<\infty\), jointly in \(y\) and \(\e\).
\end{remark}

\section{Proof for existence and uniqueness} \label{sec:EU}
\setcounter{equation}{0}
In this section, we present the proof of Theorem \ref{thm_EU}.
This section is organized as follows.
First of all, we will transform the original system \eqref{shock_0} into an alternative ODE system \eqref{sys} to which Proposition \ref{prop_FFT} can be applied.
In Section \ref{section E}, we will demonstrate that the solution to \eqref{sys} indeed gives the solution to \eqref{shock_0} to establish the existence.
The uniqueness will be addressed in Section \ref{section U}, through a direct analysis on \eqref{shock_0}.

\vspace{2mm}
In what follows, we use the following notation: \(\t_- = \t(\th_-), \m_- = \m(\th_-)\), and  \(\k_- = \k(\th_-)\) denotes the coefficients at \(\th_-\).

\vspace{2mm}
We consider the following system of ODEs: 
\begin{equation}\label{sys}
\left\{
\begin{aligned}
    w_0'(z) &= w_1(z), \\
    \e w_1'(z) &= w_2(z), \\
    \e w_2'(z) &= f(w_0(z), w_1(z), w_2(z), \e).
\end{aligned}
\right.
\end{equation}
Here \(\e\in\RR\) is a small fixed parameter and will be given by \(\e=v_+-v_-\).
However, for \eqref{sys}, the parameter \(\e\) is not necessarily assumed to be positive.
The function \(f\) is to be determined so that the system above is related to \eqref{shock_0}, and hence to \eqref{vthODE}.

To derive the system above, we first define \(g_1\) and \(g_2\) as in \eqref{gdef} through the implicit function theorem.
Since the Jacobian matrix evaluated at \((v_-, \th_-, 0, 0, 0, 0)\) is
\begin{equation}\label{hderth}
\frac{\rd (h_1, h_2)}{\rd (\thtil, \thtil')}
=
\begin{pmatrix}
    \frac{\rd h_1}{\rd \thtil} & \frac{\rd h_1}{\rd \thtil'} \\ 
    \frac{\rd h_2}{\rd \thtil} & \frac{\rd h_2}{\rd \thtil'}
\end{pmatrix}
= \begin{pmatrix}
    \frac{R}{v_-} & 0 \\ -\frac{\s_* R}{\g-1} & -\frac{\k_-}{v_-}
\end{pmatrix}
\end{equation}
and its determinant is nonzero, it follows from the implicit function theorem that there exist two smooth functions \(g_1, g_2\colon U_0\to \RR\) which satisfy \eqref{gdef}, where \(U_0\) is a neighborhood of \((v_-, 0, 0, 0)\).

We now aim to rewrite \eqref{absODE} in terms of \(w\) variables, using the relations defined in \eqref{vtilinw}:
\begin{equation}\label{finter}
\e w_2' = \Big(\frac{\rd g_1}{\rd \vtil''}\Big)^{-1}
\Big[\frac{1}{\e^2}g_2 - \frac{\rd g_1}{\rd \vtil}w_1 - \frac{\rd g_1}{\rd \vtil'}w_2\Big],
\end{equation}
provided that \(\e\neq 0\).
Notice that all \(g_i\) and their derivatives are evaluated at 
\[
(v_- + \e w_0, \e^2 w_1, \e^2 w_2, \e).
\]
Hence, we define the function \(f\) as follows:
\begin{multline} \label{f-def}
f(w_0, w_1, w_2, \e) \coloneqq
\Big(\frac{\rd g_1}{\rd \vtil''}(v_-+\e w_0, \e^2 w_1, \e^2 w_2, \e)\Big)^{-1}
\Big(\frac{1}{\e^2}g_2(v_-+\e w_0, \e^2 w_1, \e^2 w_2, \e) \\
- \frac{\rd g_1}{\rd \vtil}(v_-+\e w_0, \e^2 w_1, \e^2 w_2, \e)w_1 - \frac{\rd g_1}{\rd \vtil'}(v_-+\e w_0, \e^2 w_1, \e^2 w_2, \e)w_2\Big).
\end{multline}

The transformation to the alternative system \eqref{sys} enables us to apply Proposition \ref{prop_FFT}.
However, prior to that, it is necessary to show that the function \(f\) is not singular near \(\e=0\).

Thanks to \eqref{gdef}, it holds that \(g_1(v_-, 0, 0, 0) = \th_-\) and \(g_2(v_-, 0, 0, 0) = 0\).
Then, we need to compute all the first-order derivatives of \(g_1\) and \(g_2\) at \((v_-, 0, 0, 0)\).
Using the chain rule on \eqref{gdef}, we obtain
\begin{equation} \label{hg-ChainR}
\frac{\rd(h_1, h_2)}{\rd(\vtil, \vtil', \vtil'', \e)} +
\frac{\rd(h_1, h_2)}{\rd(\thtil, \thtil')}
\frac{\rd(g_1, g_2)}{\rd(\vtil, \vtil', \vtil'', \e)} = 0.
\end{equation}
Then, we compute all the derivatives of \(h_1\) and \(h_2\):
\begin{align*}
&\frac{\rd h_1}{\rd \vtil} 
= \s_\e^2 - \s_\e\t(\thtil)\frac{\vtil'}{\vtil^2} -\frac{R\thtil}{\vtil^2}
-\s_\e\frac{\m(\thtil)}{\vtil^2}\vtil'-2\frac{\m(\thtil)\t'(\thtil)}{\vtil^3}\vtil'\thtil' 
-2\frac{\m(\thtil)\t(\thtil)}{\vtil^3}\vtil'' 
+ 3\frac{\m(\thtil)\t(\thtil)}{\vtil^4}(\vtil')^2, \\
&\frac{\rd h_2}{\rd \vtil} 
= \s_\e\Big(\s_\e(\vtil-v_-)+\t(\thtil)\frac{\vtil'}{\vtil}\Big)\Big(\s_\e-\t(\thtil)\frac{\vtil'}{\vtil^2}\Big) 
- p_-\s_\e + p_-\t(\thtil)\frac{\vtil'}{\vtil^2} + \k(\thtil)\frac{\thtil'}{\vtil^2}, \\
&\frac{\rd h_1}{\rd \vtil'} = 
\s_\e\frac{\t(\thtil)}{\vtil} + \s_\e\frac{\m(\thtil)}{\vtil} + \frac{\m(\thtil)\t'(\thtil)}{\vtil^2}\thtil' 
-2\frac{\m(\thtil)\t(\thtil)}{\vtil^3}\vtil', 
\hspace{30mm}\frac{\rd h_1}{\rd \vtil''} = \frac{\m(\thtil)\t(\thtil)}{\vtil^2},\\
&\frac{\rd h_2}{\rd \vtil'} = \s_\e\Big(\s_\e(\vtil-v_-)+\t(\thtil)\frac{\vtil'}{\vtil}\Big)\frac{\t(\thtil)}{\vtil} 
- p_-\frac{\t(\thtil)}{\vtil},
\hspace{40mm}\frac{\rd h_2}{\rd \vtil''} = 0,\\
&\frac{\rd h_1}{\rd \e}
= \frac{\rd \s_\e}{\rd \e}\Big[2\s_\e(\vtil-v_-) + \t(\thtil)\frac{\vtil'}{\vtil} + \frac{\m(\thtil)}{\vtil}\vtil'\Big],\\
&\frac{\rd h_2}{\rd \e}
= \frac{\rd \s_\e}{\rd \e}\Big[-\frac{R}{\g-1}(\thtil-\th_-) + \frac{1}{2}\Big(\s_\e(\vtil-v_-)+\t(\thtil)\frac{\vtil'}{\vtil}\Big)^2 \\
&\qquad\qquad\qquad\qquad
+\s_\e\Big(\s_\e(\vtil-v_-)+\t(\thtil)\frac{\vtil'}{\vtil}\Big)(\vtil-v_-)- p_-(\vtil-v_-)\Big].
\end{align*}
Thus, evaluating at \((v_-, 0, 0, 0)\), we have
\begin{equation}\label{hderiv}
\frac{\rd(h_1, h_2)}{\rd(\vtil, \vtil', \vtil'', \e)}
= \begin{pmatrix}
    (\g-1)\frac{p_-}{v_-}& \s_*\frac{\t_-+\m_-}{v_-} & \frac{\m_-\t_-}{v_-^2} & 0 \\ 
    -\s_*p_-& -\frac{p_-}{v_-}\t_- & 0 & 0
\end{pmatrix},
\end{equation}
and this together with \eqref{hderth} implies that
\begin{equation}\label{gder}
\frac{\rd(g_1, g_2)}{\rd(\vtil, \vtil', \vtil'', \e)}
= \begin{pmatrix}
    -(\g-1)\frac{\th_-}{v_-} & -\frac{\s_*}{R}(\t_-+\m_-) & -\frac{\m_-\t_-}{Rv_-} & 0 \\
    0 & \frac{1}{\g-1}p_-\frac{\t_-+\g\m_-}{\k_-} & \frac{\s_*}{\g-1}\frac{\m_-\t_-}{\k_-} & 0 
\end{pmatrix}.
\end{equation}
Therefore, there exists a smaller neighborhood \(U_1\subset U_0\) of \((v_-, 0, 0, 0)\) on which \(\frac{\rd g_1}{\rd \vtil'}, \frac{\rd g_1}{\rd \vtil''}\) and \(\frac{\rd g_2}{\rd \vtil''}\) do not vanish.

Let \(M>0\) be a large constant which will be determined later. (In what follows, we will demonstrate that it depends solely on the system, especially the transport coefficients \(\t,\m,\k\) and the left-end state.)
Then, there is a small \(\e_0>0\) such that \((v_-+\e w_0, \e^2 w_1, \e^2 w_2, \e)\in U_1\) for any \((w_0, w_1, w_2)\in (-M, M)^3\) and any \(\abs{\e}< \e_0\).
Setting \(V \coloneqq (-M, M)^3 \times (-\e_0, \e_0)\), we establish that \(f\) is well-defined and smooth on \(V \setminus \{\e=0\}\).
To show that the function \(f\) can be smoothly extended to the entire \(V\), it suffices to prove that \(\frac{1}{\e^2}g_2(v_-+\e w_0, \e^2 w_1, \e^2 w_2, \e)\) can be smoothly extended to \(V\). 
To this end, we consider the Taylor expansion of \(g_2(v_-+\e w_0, \e^2 w_1, \e^2 w_2, \e)\) at \(\e=0\) with respect to \(\e\): 
\begin{multline}\label{Taylor}
g_2(v_-+\e w_0, \e^2 w_1, \e^2 w_2, \e)
= g_2(v_-, 0, 0, 0) + \e \left(w_0\frac{\rd g_2}{\rd \vtil} + \frac{\rd g_2}{\rd \e}\right)(v_-, 0, 0, 0) \\
+ \frac{\e^2}{2}\left(w_0^2\frac{\rd^2 g_2}{\rd \vtil^2} + 2w_0\frac{\rd^2 g_2}{\rd \vtil\rd \e} + 2w_1\frac{\rd g_2}{\rd \vtil'} + 2w_2\frac{\rd g_2}{\rd \vtil''} + \frac{\rd^2 g_2}{\rd \e^2}\right)(v_-, 0, 0, 0) + l.o.t..
\end{multline}
Thanks to \eqref{gder}, it holds that \(g_2(v_-, 0, 0, 0) = \frac{\rd g_2}{\rd \vtil}(v_-, 0, 0, 0) = \frac{\rd g_2}{\rd \e}(v_-, 0, 0, 0) = 0\).
Hence, we conclude that \(\frac{1}{\e^2}g_2(v_-+\e w_0, \e^2 w_1, \e^2 w_2, \e)\) is well-defined and smooth on the entire \(V\). 
Therefore, we derive the alternative (regular) system \eqref{sys} of the form \eqref{slow}.

\subsection{Proof of Theorem \ref{thm_EU}: Existence} \label{section E}
This subsection is dedicated to the existence statement of Theorem \ref{thm_EU}.
The proof consists of two steps as follows:
first, we use Proposition \ref{prop_FFT} to prove that the system \eqref{sys} admits a solution connecting the two critical points \((0, 0, 0)\) and \((1, 0, 0)\); subsequently, we demonstrate that this solution gives a desired solution to the original system \eqref{shock_0}.

\step{1} We establish that \eqref{sys} admits a solution connecting the two critical points.

\vspace{2mm}
Our starting point is to derive the critical manifold of the system \eqref{sys}; we have \(w_2=0\) and \(f(w_0, w_1, w_2, 0) = 0\).
Then, from \eqref{f-def} and \eqref{Taylor}, it is necessary to evaluate the second-order derivatives 
\[
\frac{\rd^2 g_2}{\rd \vtil^2}, \qquad \frac{\rd^2 g_2}{\rd \vtil\rd \e}, \qquad \frac{\rd^2 g_2}{\rd \e^2},
\]
at \((v_-, 0, 0, 0)\).
For this purpose, we differentiate \eqref{gdef} twice and evaluate at \((v_-, 0, 0, 0)\) as follows: using the chain rule and 
\[
\frac{\rd g_1}{\rd \e} = \frac{\rd g_2}{\rd \vtil} = \frac{\rd g_2}{\rd \e} = 0 \quad
\text{at }\, (v_-, 0, 0, 0),
\]
we find that for each \(i = 1, 2\),
\begin{equation} \label{ChainR-2}
\begin{aligned}
\frac{\rd^2 h_i}{\rd \vtil^2} + 2 \frac{\rd^2 h_i}{\rd \vtil\rd \thtil}\frac{\rd g_1}{\rd \vtil} + \frac{\rd^2 h_i}{\rd \thtil^2}\Big(\frac{\rd g_1}{\rd \vtil}\Big)^2
+ \frac{\rd h_i}{\rd \thtil} \frac{\rd^2 g_1}{\rd \vtil^2} + \frac{\rd h_i}{\rd \thtil'}\frac{\rd^2 g_2}{\rd \vtil^2} &= 0, \\
\frac{\rd^2 h_i}{\rd \e\rd \vtil} 
+ \frac{\rd^2 h_i}{\rd \thtil\rd \e}\frac{\rd g_1}{\rd \vtil}
+ \frac{\rd h_i}{\rd \thtil}\frac{\rd^2 g_1}{\rd \e\rd \vtil}
+ \frac{\rd h_i}{\rd \thtil'}\frac{\rd^2 g_2}{\rd \e\rd \vtil} &= 0, \\
\frac{\rd^2 h_i}{\rd \e^2} + \frac{\rd h_i}{\rd \thtil} \frac{\rd^2 g_1}{\rd \e^2} + \frac{\rd h_i}{\rd \thtil'}\frac{\rd^2 g_2}{\rd \e^2} &= 0.
\end{aligned}
\end{equation}
Then, it requires to calculate the following twelve second-order derivatives: 
\[
\begin{pmatrix}
\frac{\rd^2 h_1}{\rd \vtil^2} & \frac{\rd^2 h_1}{\rd \vtil\rd \thtil} & 
\frac{\rd^2 h_1}{\rd \thtil^2} & \frac{\rd^2 h_1}{\rd \e\rd \vtil} &
\frac{\rd^2 h_1}{\rd \e\rd \thtil} & \frac{\rd^2 h_1}{\rd \e^2} \\
\frac{\rd^2 h_2}{\rd \vtil^2} & \frac{\rd^2 h_2}{\rd \vtil\rd \thtil} & 
\frac{\rd^2 h_2}{\rd \thtil^2} & \frac{\rd^2 h_2}{\rd \e\rd \vtil} &
\frac{\rd^2 h_2}{\rd \e\rd \thtil} & \frac{\rd^2 h_2}{\rd \e^2}
\end{pmatrix} = 
\begin{pmatrix}
2\frac{p_-}{v_-^2} & -\frac{R}{v_-^2} & 0 & -\frac{\g(\g+1)p_-}{2v_-^2} & 0 & 0 \\
\frac{\s_*\g p_-}{v_-} & 0 & 0 & \frac{\s_*(\g+1)p_-}{4v_-} & \frac{\s_*R(\g+1)}{4v_-(\g-1)} & 0
\end{pmatrix}.
\]
This immediately implies that \(\frac{\rd^2 g_1}{\rd \e^2} = \frac{\rd^2 g_2}{\rd \e^2} = 0\), and hence, the first two equations of \eqref{ChainR-2} are equivalent to 
\[
\begin{pmatrix}
2\g\frac{p_-}{v_-^2} & -\frac{\g(\g+1)p_-}{2v_-^2} \\ \frac{\s_*\g p_-}{v_-} & 0 
\end{pmatrix}
+ \begin{pmatrix}
    \frac{R}{v_-} & 0 \\ -\frac{\s_*R}{\g-1} & -\frac{\k_-}{v_-}
\end{pmatrix}
\begin{pmatrix}
    \frac{\rd^2 g_1}{\rd \vtil^2} & \frac{\rd^2 g_1}{\rd \vtil\rd \e} \\
    \frac{\rd^2 g_2}{\rd \vtil^2} & \frac{\rd^2 g_2}{\rd \vtil\rd \e}
\end{pmatrix} = 0.
\]
Then, an elementary calculation yields that
\[
\begin{pmatrix}
    \frac{\rd^2 g_1}{\rd \vtil^2} & \frac{\rd^2 g_1}{\rd \vtil\rd \e} \\
    \frac{\rd^2 g_2}{\rd \vtil^2} & \frac{\rd^2 g_2}{\rd \vtil\rd \e}
\end{pmatrix}
= \frac{v_-^2}{R\k_-}\begin{pmatrix}
-\frac{k_-}{v_-} & 0 \\ \s_*\frac{R}{\g-1} & \frac{R}{v_-}
\end{pmatrix}
\begin{pmatrix}
2\g\frac{p_-}{v_-^2} & -\frac{(\g+1)\s_*^2}{2v_-} \\ \frac{\s_*\g p_-}{v_-} & 0 
\end{pmatrix}
= \begin{pmatrix}
    -\frac{2\g p_-}{Rv_-} & \frac{(\g+1)\s_*^2}{2R} \\ 
    \frac{\g+1}{\g-1}\frac{\s_*\g p_-}{\k_-} & -\frac{\g+1}{\g-1}\frac{\s_*\g p_-}{2\k_-}
\end{pmatrix}.
\]
Note that all the derivatives above were evaluated at \((v_-, 0, 0, 0)\).
The critical manifold of the system \eqref{sys}, which is given by \(w_2 = 0\) and \(f(w_0, w_1, 0, 0) = 0\), can be represented by
\begin{align*}
0 &= \frac{w_0^2}{2}\frac{\rd^2 g_2}{\rd \vtil^2} + w_0 \frac{\rd^2 g_2}{\rd \vtil\rd \e} + w_1\frac{\rd g_2}{\rd \vtil'} - w_1\frac{\rd g_1}{\rd \vtil} \\
&= \frac{\g+1}{\g-1}\frac{\s_*\g p_-}{2\k_-}(w_0^2-w_0) 
+ \Big(\frac{1}{\g-1}\frac{\t_-+\g\m_-}{\k_-}p_- + \frac{\g-1}{R}p_-\Big)w_1.
\end{align*}
Therefore, the critical manifold \(M_0\) is a one-dimensional manifold given by the graph
\begin{equation}\label{Critmfd}
M_0 = \left\{(w_0, w_1, w_2)\in V \middle| w_2 = 0, w_1 = A(w_0-w_0^2)\right\},
\end{equation}
where \(A\) is the constant defined in \eqref{Adef}.

\vspace{2mm}
We may now proceed to fix the size of the domain \(V\).
It is required that \(M > \max\{1,A\}\) so that the set \(M_0\cap \{w_0\in[-1/2,3/2]\}\) is contained in the projection of \(V\) onto \(\{\e=0\}\).

\vspace{2mm}
To apply Proposition \ref{prop_FFT}, it is necessary to verify that the critical manifold \(M_0\) is normally hyperbolic.
To this end, we claim that each entry on the second row of
\begin{align*}
\left.\begin{pmatrix}
    \frac{\rd w_2}{\rd w_1} & \frac{\rd w_2}{\rd w_2} \\ 
    \frac{\rd}{\rd w_1}f(w_0, w_1, w_2, \e) & \frac{\rd}{\rd w_2}f(w_0, w_1, w_2, \e)
\end{pmatrix}\right|_{\e = 0}
\end{align*}
is negative.
Considering \eqref{f-def} and \eqref{Critmfd}, when the derivative is applied to \(\frac{\rd g_1}{\rd \vtil''}\), the remaining factor \((\frac{1}{\e^2}g_2 - \frac{\rd g_1}{\rd \vtil}w_1 - \frac{\rd g_1}{\rd \vtil'}w_2)\) vanishes on \(M_0\).
Notice that \(\frac{\rd g_1}{\rd \vtil''}\) is always negative on \(M_0\) due to the choice of \(U_1\).
For the \(w_1\)-derivative, we examine \(\frac{\rd g_2}{\rd \vtil'}\) and \(-\frac{\rd g_1}{\rd \vtil}\): since both are positive, the \(w_1\)-derivative is negative.
For the \(w_2\)-derivative, we analyze \(\frac{\rd g_2}{\rd \vtil''}\) and \(-\frac{\rd g_1}{\rd \vtil'}\): since both of them are positive, the \(w_2\)-derivative is also negative.
Therefore, the trace of the matrix above is negative and its determinant is positive.
This excludes the case where any eigenvalue lies on the imaginary axis.

\vspace{2mm}
We now confirm the existence of two isolated critical points of \eqref{sys}, namely, \((0, 0, 0)\) and \((1, 0, 0)\), for any \(\e\neq0\) as follows.
Note that \((w_0, w_1, w_2)\) is a critical point of \eqref{sys} if and only if \(w_1 = w_2 = 0\) and \(g_2(v_-+\e w_0, 0, 0, \e) = 0\).
Then, since
\[
h_i(v_\pm, \th_\pm, 0, 0, 0, \e) = 0 ~\text{ for }~ i = 1, 2,
\]
it holds that \((0, 0, 0)\) and \((1, 0, 0)\) are two critical points.
Further, if \((w_0, 0, 0)\) is a critical point, then it follows that for \(\th = g_1(v_-+\e w_0, 0, 0, \e)\),
\[
h_i(v_-+\e w_0, \th, 0, 0, 0, \e) = 0 \quad \text{for }\, i = 1, 2.
\]
This implies the following:
\begin{align*}
\e\s_\e^2 w_0 + \frac{R\th}{v_-+\e w_0} - p_- = 0, \qquad
-\frac{\s_\e R}{\g-1}(\th-\th_-) + \e^2\frac{\s_\e^3}{2}w_0^2 - \e\s_\e p_-w_0 = 0.
\end{align*}
Then, eliminating \(\th\), we obtain from \eqref{seval} that
\begin{align*}
0 &= 
\e^2\frac{\s_\e^3}{2}w_0^2 - \e\s_\e p_-w_0 + \frac{\s_\e R}{\g-1}\th_-
- \frac{\s_\e}{\g-1}(v_-+\e w_0)(p_--\e\s_\e^2w_0) \\
&= -\e^2\s_\e^3\frac{\g+1}{2(\g-1)}w_0(1-w_0),
\end{align*}
which proves that \(w_0 = 0\) or \(1\).
Thus, for each nonzero \(\e\), there are only two critical points \((0,0,0)\) and \((1,0,0)\) in the domain \((-M, M)^3\).

\vspace{2mm}
The remaining argument is the same as \cite{EEKO}.
Firstly, we consider a closed interval \(K\) such that \([0,1] \subset K \subset [-1/2,3/2]\).
Then, we apply Proposition \ref{prop_FFT} and Remark \ref{rmk-FFT} to obtain a locally invariant manifold \(M_\e\) as a graph of a \(\Ccal^3\) function of \(w_0\) for each (small) nonzero \(\e\).
More explicitly, we choose sufficiently regular (jointly in \(w_0\) and \(\e\)) functions \(s_1\) and \(s_2\) of two variables \(w_0\in K\) and \(\e\) such that
\begin{equation}\label{invarmfd}
\{(w_0, w_1, w_2)|w_0\in K, w_1 = Aw_0(1-w_0)+\e s_1(w_0, \e), w_2 = \e s_2(w_0, \e)\}.
\end{equation}
This is a parametrization of the invariant manifold \(M_\e\) for \eqref{sys} by a small parameter \(\e\neq0\).
Thereafter, using classical ODE theory for the autonomous case and the Cauchy-Lipschitz theorem, we demonstrate that \eqref{sys} admits a solution \((w_0, w_1, w_2)\colon \RR\to (-M, M)^3\) for each \(\e\neq0\), which connects the two critical points \((0, 0, 0)\) and \((1, 0, 0)\). 
It is worth noting that these two points should be on the invariant manifold, i.e., 
\begin{equation}\label{szeros}
s_1(0, \e) = s_1(1, \e) = s_2(0, \e) = s_2(1, \e) = 0.
\end{equation}

\step{2} We construct a solution to the original system \eqref{shock_0} from \((w_0, w_1, w_2)\) which solves alternative system \eqref{sys}.

\vspace{2mm}
In accordance with \eqref{vtilinw}, we define 
\begin{equation}\label{solback}
\begin{aligned}
\vtil(\x) &\coloneqq \e w_0(\e\x) + v_-, \qquad
\thtil(\x) \coloneqq g_1(\vtil(\x), \vtil'(\x), \vtil''(\x), \e), \\
\util(\x) &\coloneqq u_- - \s_\e(\vtil(\x)-v_-)-\t(\thtil(\x))\frac{\vtil'(\x)}{\vtil(\x)}.
\end{aligned}
\end{equation}
It is obvious that these functions satisfy \eqref{uexplicit}.
Moreover, using \eqref{finter} we obtain
\begin{equation}\label{solback2}
\thtil'(\x)
= \frac{\rd g_1}{\rd \vtil}\e^2w_1(\e\x) + \frac{\rd g_1}{\rd \vtil'}\e^2w_2(\e\x)
+ \frac{\rd g_1}{\rd \vtil''}\e^3w_2'(\e\x)
= g_2(\vtil(\x), \vtil'(\x), \vtil''(\x), \e),
\end{equation}
and hence \eqref{vthODE} holds.
Therefore, we conclude that \((\vtil, \util, \thtil)\) in \eqref{solback} is indeed a solution to the original system \eqref{shock_0} with \(\s = \s_\e\), which connects the two end states \((v_-, u_-, \th_-)\) and \((v_+, u_+, \th_+)\).

\vspace{2mm}
Note that since \(w_0\) is increasing, the monotonicity of \(\vtil\) follows from its definition \eqref{solback}.
The monotonicity of \(\util\) and \(\thtil\) is a consequence of \eqref{ratio-vu} and \eqref{ratio-vth}, whose proofs are deferred to Section \ref{sec:property}.

\subsection{Proof of Theorem \ref{thm_EU}: Uniqueness} \label{section U}
In Section \ref{section E}, we established the existence of a solution to \eqref{intODE} which connects \((v_-, u_-, \th_-)\) to \((v_+, u_+, \th_+)\).
This subsection is devoted to proving the uniqueness of the solution derived in Section \ref{section E}.

\vspace{2mm}
To show the uniqueness, we follow the argument developed in the constant dissipation coefficient case \cite{EEKO}, which is based on the following proposition of the unstable manifold at a critical point.

\begin{proposition} \cite[Theorems 9.4, 9.5]{teschl2012ordinary}  \label{TeschlUnique}
Consider a system of ODEs \(z' = F(z)\), where  \(F \in \Ccal^k(\RR^n ; \RR^n)\) and $x_0$ is a critical point, i.e., \(F(x_0) = 0\). Let \(W^-(x_0)\) be the unstable manifold composed of all points converging to \(x_0\) for \(t\to -\infty\). 
Let \(E^+\) (resp. \(E^-\)) denote the linear subspace spanned by eigenvectors of \(D_{x_0} F\) corresponding to negative (resp. positive) eigenvalues. 
Suppose \(D_{x_0} F \) has no eigenvalues on the imaginary axis. \\
Then, there are a neighborhood \(U(x_0) = x_0 + U\) of \(x_0\) and a function \(h^-\in \Ccal^k(E^-\cap U; E^+)\) such that 
\[
W^-(x_0)\cap U(x_0) = \{x_0+a+h^-(a)\mid a\in E^-\cap U\}.
\]
Here, both \(h^-\) and its Jacobian matrix vanish at \(0\).
\end{proposition}
This proposition says that \(W^-(x_0)\) is a \(\Ccal^k\)-manifold whose dimension is equal to \(\dim E^-\) and is tangent to the affine space \(x_0 + E^-\) at \(x_0\).

\vspace{2mm}
The system \eqref{intODE} can be written in the following form: 
\begin{equation} \label{intODE2}
\left\{
\begin{aligned}
    & \vtil' = \frac{\vtil}{\t(\thtil)}(-\s_\e(\vtil-v_-)-(\util-u_-)), \\
    & \util' = \frac{\vtil}{\m(\thtil)}(-\s_\e(\util-u_-)+(\ptil-p_-)), \\
    & \thtil' = \frac{\vtil}{\k(\thtil)}\Big[-\s_\e\Big(\frac{R}{\g-1}(\thtil-\th_-)+\frac{1}{2}(\util-u_-)^2\Big)+p_-(\util-u_-)\Big].
\end{aligned}
\right.
\end{equation}
Thanks to \eqref{end-con}, it is obvious that the two end states \((v_-,u_-,\th_-)\) and \((v_+,u_+,\th_+)\) are the critical points of the system \eqref{intODE2}.
We proceed to examine the sign of eigenvalues of the linearized system at the critical point \((v_-,u_-,\th_-)\) so as to determine the dimension of the unstable manifold.
The Jacobian matrix at \((v_-, u_-, \th_-)\) is given by
\begin{equation} \label{Jacobian}
J(\vtil, \util, \thtil)|_{(v_-, u_-, \th_-)}
= \begin{pmatrix}
    -\frac{\s_\e v_-}{\t_-} & -\frac{v_-}{\t_-} & 0 \\
    -\frac{p_-}{\m_-} & -\frac{\s_\e v_-}{\m_-} & \frac{R}{\m_-} \\
    0 & \frac{p_-v_-}{\k_-} & -\frac{\s_\e R}{\g-1}\frac{v_-}{\k_-}
\end{pmatrix}.
\end{equation}
Then, the determinant and trace of the matrix \eqref{Jacobian} are as follows:
\begin{align*}
\det J((v_-, u_-, \th_-)) &= \frac{1}{\t_-\m_-\k_-}\frac{R\s_\e v_-^2}{\g-1}(\g p_--\s_\e^2v_-) > 0, \\
\tr J((v_-, u_-, \th_-)) &= -\s_\e v_-\Big(\frac{R}{\g-1}\frac{1}{\k_-} + \frac{1}{\m_-} + \frac{1}{\t_-}\Big)<0.
\end{align*}
We used \eqref{seval} for the first inequality.
Given that the determinant is positive and the trace is negative, it follows that it has one positive real eigenvalue and two eigenvalues with negative real part.
Thus, the dimension of the unstable manifold is \(1\).

\vspace{2mm}
We now show the desired uniqueness as follows:
as established earlier, there exists a monotone solution \((\vtil, \util, \thtil)\) to \eqref{intODE}, connecting the two states \((v_-, u_-, \th_-)\) to \((v_+, u_+, \th_+)\).
This solution trajectory lies locally on the one-dimensional unstable manifold \(W^-((v_-, u_-, \th_-))\).
More precisely, the intersection between \(W^-((v_-, u_-, \th_-))\) and \(\{v\ge v_-, u \le u_-, \th\le \th_-\}\) coincides with the trajectory \((\vtil, \util, \thtil)\) in a neighborhood of \((v_-, u_-, \th_-)\).
Furthermore, any other solution that connects \((v_-, u_-, \th_-)\) and \((v_+, u_+, \th_+)\) must lie on the unstable manifold \(W^-((v_-, u_-, \th_-))\).
In particular, if such a solution is monotone, it must intersect the trajectory \((\vtil, \util, \thtil)\).
Therefore, by the uniqueness theorem for Lipschitz continuous autonomous system of ODE, we conclude that \((\vtil, \util, \thtil)\) is the unique solution in the class of monotone solutions, up to spatial translation.

\section{Proof of Theorem \ref{thm_estimates}}
\label{sec:property}
\setcounter{equation}{0}

Now we present the proof of Theorem \ref{thm_estimates}.
Most of the estimates proceed in a manner similar to \cite{EEKO}, as the invariant manifold---hence the functions \(s_1\) and \(s_2\)---exhibits the same regularity as in the constant dissipation coefficient case.
The primary distinction lies in the fact that we now aim to determine the precise ratio between \(\vt\) and \((v_+-\vt)(\vt-v_-)\).

\vspace{2mm}
\(\bullet\) \textbf{Notation:} In what follows, \(C\) denotes a positive constant that may change from line to line, yet remains independent of \(\e\).

\vspace{2mm}
To prove Theorem \ref{thm_estimates}, the argument will be developed in the following sequence.
First of all, we obtain \((\abs{\vt''(\x)},\abs{\ut''(\x)},|\tht''(\x)|) \le C\e\abs{\vt'(\x)}\)  based on the fact that the solution trajectory lies on the invariant manifold, whose shape is known.
We then establish \eqref{ratio-vu} and \eqref{ratio-vth}, from which \eqref{derivsign} and \eqref{2bound1} follow immediately.
Subsequently, we prove \eqref{tail}, which, together with Gronwall's inequality, yields \eqref{decay} and \eqref{derivdecay}.

\vspace{2mm}
First of all, we recall \eqref{sys} and \eqref{vtilinw}, and thus we have
\begin{align} 
\vt'(\x) &= \e^2 w'_0(z),
&\vt''(\x) &= \e^2 w_2(z) = \e^3 w_1'(z), \label{vw12}\\
\vt'''(\x) &= \e^4w''_1(z),
&\vt''''(\x) &= \e^5w'''_1(z). \label{vw34}
\end{align}
Moreover, it follows from \eqref{invarmfd} that
\begin{align} \label{w1=}
w_1 = A(w_0-w_0^2) + \e s_1(w_0, \e), \quad w_0\in K,
\end{align}
for a compact set \(K\) and a smooth function \(s_1\), which directly implies \(\abs{w_1} = \abs{w'_0} \le C\).\\
Hence, \eqref{vw12} simply yields that \(\abs{\vt'(\x)}\le C\e^2\).
We continue by differentiating \eqref{w1=}:
\[
w'_1 = A(1-2w_0)w'_0 + \e D_1 s_1(w_0, \e)w'_0.
\]
Then, since \(D_1 s_1\) is continuous and \(w_0\) is confined to \(K\), it follows that \(\abs{w'_1(z)}\le C \abs{w'_0(z)}\), which together with \eqref{vw12} gives that \(\abs{\vt''(\x)}\le C\e\abs{\vt'(\x)}\).

It also requires to obtain the bounds for \(\abs{\ut''}\) and \(|\tht''|\).
Due to \eqref{solback}, we need to control the higher-order derivatives of \(\vt\).
We differentiate \eqref{w1=} once and twice more to observe
\begin{align*}
w''_1 &= -2Aw_1w_0' + A(1-2w_0)w_1' + \e D_{11}s_1(w_0,\e)w_1w_0' + \e D_1s_1(w_0,\e)w_1',\\
w'''_1 &= -6A w_1w_1' + A(1-2w_0)w_1'' + \e D_{111}s_1(w_0,\e)w_1^2w_0'\\
&\qquad+3\e D_{11}s_1(w_0,\e)w_1w_1'
+\e D_1s_1(w_0,\e)w_1''.
\end{align*}
Then, since derivatives of \(s_1\) up to third order are continuous and \(w_0\) is confined to \(K\) either, it holds that \(\abs{w''_1(z)}\le C\abs{w'_0(z)}\) and \(\abs{w'''_1(z)}\le C\abs{w'_0(z)}\).
This together with \eqref{vw34} implies \(\abs{\vt'''(\x)}\le C\e^2\abs{\vt'(\x)}\) and \(\abs{\vt''''(\x)}\le C\e^3\abs{\vt'(\x)}\).
Now, differentiating \eqref{solback2}, we have
\begin{align*}
\tht''(\x) = \frac{\rd g_2}{\rd \vt}\vt' + \frac{\rd g_2}{\rd \vt'}\vt'' + \frac{\rd g_2}{\rd \vt''}\vt'''.
\end{align*}
In view of the Taylor expansion on \(\frac{\rd g_2}{\rd \vt}(v_-+\e w_0, \e^2 w_1, \e^2 w_2, \e)\) with respect to \(\e\), and noting by \eqref{gder} that \(\frac{\rd g_2}{\rd \vt}\) vanishes at \((v_-,0,0,0)\), it follows that \(\frac{\rd g_2}{\rd \vt}\) is of order \(\Ocal(\e)\).
Similarly, \(\frac{\rd g_2}{\rd \vt'}\) and \(\frac{\rd g_2}{\rd \vt''}\) are of order \(\Ocal(1)\).
Thus, we obtain that \(|\tht''(\x)|\le C\e\abs{\vt'(\x)}\).\\
The situation is similar for \(\ut''\).
We differentiate \eqref{solback} twice to find that
\begin{align*}
\ut'' &=
-\s_\e \vt''
-\t''(\tht)(\tht')^2 \frac{\vt'}{\vt}
-\t'(\tht)\tht'' \frac{\vt'}{\vt}
-2\t'(\tht)\tht'\frac{\vt''\vt-(\vt')^2}{\vt^2}\\
&\qquad
-\t(\tht)\frac{\vt'''(\vt)^2 -3\vt''\vt'\vt+2(\vt')^3}{(\vt)^3}
\end{align*}
which immediately implies \(\abs{\ut''(\x)}\le C\e\abs{\vt'(\x)}\).

\vspace{2mm}
Now we show \eqref{ratio-vu} and \eqref{ratio-vth} as follows.
Using \eqref{shock_0}\(_1\), it follows that
\[
\abs{\ut' + \s_\e\vt'} = \abs{\t'(\tht)\tht'\frac{\vt'}{\vt}+\t(\tht) \frac{\vt''}{\vt}-\t(\tht)\frac{(\vt')^2}{\vt^2}} \le C\e\abs{\vt'}.
\]
This together with \eqref{sm1} implies \eqref{ratio-vu}.
Furthermore, for \eqref{ratio-vth}, we utilize \eqref{gder} to observe
\[
\abs{\tht' + (\g-1)\frac{\th_-}{v_-}\vt'} \le C\e \abs{\vt'}.
\]
Note that these two estimates with \(\vt'>0\) imply that \(\abs{\vt'}\sim \abs{\ut'}\sim |\tht'|\) and \eqref{derivsign}.

\vspace{2mm}
We turn to the proof of \eqref{tail}.
Recall the fact that the two critical points \((0,0,0)\) and \((1,0,0)\) lie on the invariant manifold and thus \eqref{szeros} holds.
Since the two functions \(s_1\) and \(s_2\) are sufficiently regular, we choose two continuous functions \(\stil_1\) and \(\stil_2\) such that
\[
s_1(w_0, \e) = w_0(1-w_0)\stil_1(w_0, \e), \quad 
s_2(w_0, \e) = w_0(1-w_0)\stil_2(w_0, \e).
\]
Then, using \eqref{vtilinw} and \eqref{invarmfd}, we obtain
\[
\frac{\vt'}{(v_+-\vt)(\vt-v_-)} = \frac{\e^2 w'_0(z)}{\e^2 w_0(z)(1-w_0(z))} = \frac{w_1(z)}{w_0(z)(1-w_0(z))}
= A + \e\stil_1(w_0, \e).
\]
This establishes that 
\[
\abs{\frac{\vt'}{(v_+-\vt)(\vt-v_-)} - A} \le C\e,
\]
from which the desired conclusion \eqref{tail} follows.
Moreover, building up on \eqref{tail}, the standard Gronwall's lemma argument (for example, see \cite[Lemma 2.1]{KV21}) implies the following:
\begin{equation}
\begin{aligned}\label{vdecay}
\abs{\vt(\x)-v_-}\le C\e e^{-C\e\abs{\x}} \,\text{ for }\, \x \le 0, \qquad
\abs{\vt(\x)-v_+}\le C\e e^{-C\e\abs{\x}} \,\text{ for }\, \x \ge 0.
\end{aligned}
\end{equation}
Then, plugging these estimates into \eqref{tail}, we obtain
\[
\abs{\vt'(\x)} \le C\e^2e^{-C\e\abs{\x}},
\]
and this together with \eqref{ratio-vu} and \eqref{ratio-vth} gives the full strength of \eqref{derivdecay}.\\
In addition, applying \eqref{vdecay} and \eqref{derivdecay} to \eqref{intODE}\({}_1\) and \eqref{intODE}\({}_3\), we obtain that 
\begin{align*}
|\ut(\x)-u_-|&\le C\e e^{-C\e\abs{\x}} \,\text{ for }\, \x \le 0, 
&|\ut(\x)-u_+|&\le C\e e^{-C\e\abs{\x}} \,\text{ for }\, \x \ge 0, \\
|\tht(\x)-\th_-|&\le C\e e^{-C\e\abs{\x}} \,\text{ for }\, \x \le 0, 
&|\tht(\x)-\th_+|&\le C\e e^{-C\e\abs{\x}} \,\text{ for }\, \x \ge 0.
\end{align*}
This completes the proof of \eqref{decay} and Theorem \ref{thm_estimates}. \qed

\vspace{2mm}
Corollary \ref{Jaccor} can be deduced from the above estimates as follows.
Since \(y(\x)=\frac{\vt(\x)-v_-}{v_+-v_-}=\frac{\vt(\x)-v_-}{\e}\), it simply holds that
\[
\frac{1}{y(1-y)}\frac{dy}{d\x}
= \frac{\e}{(v_+-\vt)(\vt-v_-)}\vt'.
\]
Therefore, the preceding estimate \eqref{tail} directly implies \eqref{Jacineq}:
\[
\abs{\frac{1}{y(1-y)}\frac{dy}{d\x} - A\e}
= \e\abs{\frac{\vt'-A(v_+-\vt)(\vt-v_-)}{(v_+-\vt)(\vt-v_-)}}
\le C\e^2.
\]
This completes the proof of Corollary \ref{Jaccor}. \qed

\subsection*{Conflict of interest}
We have no conflict of interest to declare to this work.

\subsection*{Data availability statement}
We do not analyze or generate any datasets, because our work proceeds within a theoretical and mathematical approach.

\bibliographystyle{plain}
\bibliography{reference}

\end{document}